\documentclass[leqno,12pt]{article}
\usepackage{latexsym}
\usepackage{hyperref,amsmath,calligra,yfonts}
\usepackage{amssymb}
\usepackage{bbm}
\usepackage{color}
\usepackage[latin1]{inputenc}

\usepackage{graphicx}
\usepackage{epsfig}

 0

\newcommand{\R}{\mathbb{R}}

\def \be{\begin{eqnarray*}}
\def \ee{\end{eqnarray*}}
\def \ben{\begin{eqnarray}}
\def \een{\end{eqnarray}}

\newcommand{\B}{\mathcal{B}}

\newcommand{\tr}{\mathrm{tr}}

\newcommand{\br}{\mathcal{B}}
\newcommand{\usn}{a_{n}}
\newcommand{\Int}{\mathrm{int}\,}
\newcommand{\Clo}{\mathrm{clo}\,}
\newtheorem{lem}{Lemma}[section]
\newtheorem{defi}[lem]{Definition}
\newtheorem{theo}[lem]{Theorem}
\newtheorem{cor}[lem]{Corollary}
\newtheorem{prop}[lem]{Proposition}

\newtheorem{rem}[lem]{Remark}

\newcommand{\schw}{\stackrel{\mathcal{D}}{\longrightarrow}}

\makeatletter\@addtoreset{equation}{section}\makeatother

\newcommand{\ba}{\begin{array}}
\newcommand{\ea}{\end{array}}
\newcommand{\beqohne}{\begin{eqnarray*}}
\newcommand{\eeqohne}{\end{eqnarray*}}
\newcommand{\beohne}{\begin{equation*}}
\newcommand{\eeohne}{\end{equation*}}

\def\3{\ss}

\def \w{{\tt w}}
\def \u{{\tt u}}
\def \v{{\tt v}}
\def \h{{\tt h}}

\newcommand{\beq}{\begin{equation}}
\newcommand{\eeq}{\end{equation}}

\def\QED{\hfill\vrule height 1.5ex width 1.4ex depth -.1ex \vskip20pt}
\def\sn{^{(N)}}

\newcommand{\lal}{\langle\langle}
\newcommand{\rar}{\rangle\rangle}

\newcommand{\spann}{ \mbox{\sl span} \ }

\def \sur#1#2{\mathrel{\mathop{\kern 0pt#1}\limits^{#2}}}

\DeclareMathAlphabet{\mathcalligra}{T1}{calligra}{m}{n}

\begin{document}

\title{Operator-valued spectral  measures and large deviations}

\author{
{\small Fabrice Gamboa}\\
{\small Universit\'e Paul Sabatier}\\
{\small Institut de Math\'ematiques de Toulouse}\\
{\small 118 route de Narbonne}\\
{\small 31062 Toulouse Cedex 9, France}\\
{\small e-mail: gamboa@math.univ-toulouse.fr}\\ 
\and{\small Alain Rouault}\\
{\small Universit\'e Versailles-Saint-Quentin}\\
{\small LMV UMR 8100}\\
{\small 45 Avenue des Etats-Unis}\\
{\small 78035-Versailles Cedex France}\\
{\small e-mail: alain.rouault@uvsq.fr}
}

\maketitle

\begin{abstract}
Let
$\mathfrak H$  be a  
Hilbert space, let $U$ be a unitary operator on $\mathfrak H$ and let $\mathfrak K$ be a cyclic subspace for $U$.
 The spectral measure of the pair $(U, \mathfrak K)$ is an 
operator-valued measure $\mu_{\mathfrak K}$ on the unit circle $\mathbb T$ such that
\[\int_\mathbb T z^k d\mu_{\mathfrak K}(z)
 = \left( P_\mathfrak K U^k\right)_{\!\upharpoonright \mathfrak K}\ , \ \forall \ k\geq 0 \,\]
where $P_{\mathfrak K}$ and $\upharpoonright\!\mathfrak K$ are the projection  and restriction on $\mathfrak K$, respectively. 
When $\mathfrak K$ is one dimensional, $\mu$ is a scalar probability measure. In this case, 
if $U$ is picked at random from the unitary group $\mathbb U(N)$ 
under the Haar measure, then any fixed $\mathfrak K$ is almost surely cyclic for $U$. 
Let $\mu\sn$ be the  random spectral (scalar) measure
 of $(U, \mathfrak K)$. 
The sequence  $(\mu\sn)$ was studied  
 extensively, in the regime of large $N$. It converges to the Haar measure $\lambda$ on $\mathbb T$ and satisfies the Large Deviation Principle at scale $N$ with a good rate function which is  the reverse Kullback information with respect to $\lambda$  
(\cite{gamboacanonical}). The purpose of the present paper is to give an extension of this result for general $\mathfrak K$ (of fixed finite dimension $p$) and eventually to offer a projective statement (all $p$ simultaneously), at the level of 
operator-valued spectral measures in infinite dimensional spaces.
 \end{abstract}

\medskip
\section{Introduction}
\subsection{The influence of the mathematical  work of Studden on our research}
A significant part of the mathematical contribution of W.J. Studden  relies on moment problems or more generally on generalized moment problems for $T$-{\it systems}. The first author of the present paper first met  the $T$-{\it systems} during his Ph.D preparation by the 
fascinating reading of  two books on moment problems. The first one is the book of Krein and Nudel'man \cite{KreNud} dealing mainly with the Markov moment problem. The second one is the book of Karlin and Studden \cite{KarStu} that offers a beautiful journey inside  the continent of $T$-{\it systems} properties. The reading of these two books has whetted our interest for the literature on moment problems and by the middle of the nineties  we came across a very interesting paper of Chang, Kempermann and Studden \cite{Changetal} on the asymptotic behaviour of randomized moment sequences. This seminal paper  gives a very nice Borel Poincar\'e like theorem for moment sequences of  probability measures on the unit interval and has been quite motivating for at least the ten last years of our researches. The probabilized  moment space frame developed therein led to many papers written by many authors (see for example \cite{Gup}, \cite{Gamboa}, \cite{LozEJP}, \cite{DeGa},  \cite{denag09}, \cite{118218FGAR} ). One of the main ingredient tool for the study of probabilized moment spaces is the parametrization of these spaces by the canonical moments. Roughly speaking, under natural probability measures these parameters
have a joint product law with beta marginal distributions. They are also very intriguing nice mathematical objects with a lot of properties that we learned from the excellent book of Dette and Studden \cite{DeSt97}. Moreover, studying the exhaustive books of Simon (\cite{simon05}), we realized that canonical moments, also called Verblunsky coefficients,  are quite important objects in complex analysis and spectral theory. At this time the second author of the present paper was working on the asymptotic properties of the determinant of classical random matrix ensembles (\cite{Al1}). 
Surprisingly,  by the Bartlett formula, the distribution of these random determinants involves product of independent beta  random variables having similar parameters as those found in the randomized moment problem.  Observing this analogy,  we discovered a connection between random moments  and  spectral measures of classical random matrix ensembles (\cite{gamboacanonical}, \cite{FGAR}). The present paper is a matricial extension of the  
asymptotic studies conducted in the latter papers 
  dealing with scalar random spectral measures.  This extension has been possible thanks to two significant contributions of Dette and Studden in the field of matricial moment problems (see \cite{destu02}, and \cite{destu03}).

We never had the opportunity to meet W.J. Studden but we wish to pay here a tribute to this creative mathematician that had often enlighten the  paths  of our researches. 
\subsection{Introduction to this paper} 
To capture the asymptotic behavior of large dimensional unitary random matrices, the usual statistic is the empirical spectral distribution, providing equal weights to all eigenvalues
\[\mu_\u = \frac{1}{N}\sum_{k=1}^N \delta_{\lambda_k}\,.\]  More recently, some authors used another random probability measure based on eigenvalues and eigenvectors (\cite{Killip1}, \cite{bai1}, \cite{gamboacanonical}). If $U$ is a  unitary matrix and $e$ is  a fixed vector (assumed to be cyclic),  the so-called spectral measure $\mu_{\w, 1}$ of the pair $(U,e)$ may be defined through its algebraic moments. Indeed, for all $n\in\mathbb Z$
\[\langle e , U^ne\rangle = \int_{\mathbb T} z^n d\mu_{\w, 1} (z)\,.\]
Here,  $\mathbb T$ is  the unit circle $\{ z \in \mathbb C : |z| = 1\}$.
The measure $\mu_{\w, 1}$ is finitely supported on the eigenvalues of $U$, we may write
\[\mu_{\w , 1} = \sum_{k=1}^N \w_k\delta_{\lambda_k}\]
 where  $\w_k$ is the square of the scalar product  of $e$ with a unitary  eigenvector associated with $\lambda_k$.  The latter measure carries more information than the former. The weights in $\mu_{\w , 1}$ are blurred footprints of the eigenvectors of $U$. To make these footprints unblurred,
it is then tempting to try to increase the {\it dimension} by {\it projecting} $U$  on a fixed subspace of dimension $p$  instead on the span generated by the single vector $e$.  We obtain a matrix  measure.  This is what we will do in this paper.   
 Actually we may even go back to the representation given by the spectral theorem  (see \cite{DuSch})
\[U = \int_{\mathbb T} z E_U(dz)\]
where  $E_U$ is the spectral measure of $U$ (or resolution of the identity for $U$). 
 In our work, we sample $U$ according to the Haar distribution on $\mathbb U (N)$ and  we study the random object $E_U$. 
 
 The paper is organized as follows. In the next section, we first frame our paper by giving the main notations and definitions needed further. Then, we recall some facts on unitary matrices and matrix orthogonal polynomials on the unit circle. We also show technical results on these objects that will be useful later.   In Section \ref{srando} we first study the effects of the randomization of the unitary matrices on the object defined in Section \ref{spre}. In particular, our approach merely simplifies the proof of the asymptotic normality for a fixed corner extracted in the unitary ensemble given in \cite{krishna09}. Further, we show large deviation theorems both for matrix random spectral measure and their infinite dimensional lifting. All proofs are postponed to Section \ref{sproofs}.
\section{Preliminaries}
\label{spre}
\subsection{Some notations and definitions}
\label{susono}
To begin with, let us give some definitions and notations.
 Let $\mathbb N = \{1,2, \dots\}$ and 
 ${\mathcal H}=\ell^2_{\mathbb C}(\mathbb N)$.  For $i\geq 1$, let 
$e_i = (0, \dots, 0, 1, 0, \dots)$ be the $i$-th element of the canonical basis of  $\mathcal H$
 and for $p \geq 1$ let 
 ${\mathcal H}_p$ be 
    $\spann \{e_1, \cdots, e_p\}$.  
We define several sets of matrices with complex entries:
\begin{itemize}   
\item ${\mathbb M}_{p,n}$, the set of $p\times n$ matrices with complex entries,
%
%
\item $\mathbb U (n)$, the set of $n\times n$ unitary matrices,
\end{itemize}
At last,  $I_p$ denotes the $p\times p$ identity  matrix on $\mathcal H_p$.

 \subsection{Operator-valued and spectral measures}
 \label{suopsp}
Let  $\mathcal B(\mathbb T)$ denote the Borel $\sigma$-algebra on $\mathbb T$.  Let $\mathfrak H$ be  a separable Hilbert space, $I_{\mathfrak H}$ be the identity in $\mathfrak H$ and 
 $H(\mathfrak H)$ 
 be the algebra of bounded Hermitian endomorphisms on  $\mathfrak H$. 
\begin{defi}
$\;$
\begin{enumerate}
\item
A mapping $\Sigma : \mathcal B(\mathbb T) \rightarrow H(\mathfrak H)$ is called an operator-valued measure in $\mathfrak H$ if 
\begin{enumerate}
\item $\Sigma(\emptyset) = 0$;
\item $\Sigma$ is non-negative i.e. if $\Sigma(\Delta) \geq 0$ for $\Delta \in \mathcal B(\mathbb T)$;
\item $\Sigma$ is strongly countably additive, i.e., if $\Delta = \cup_{j=1}^\infty \Delta_j$ is a disjoint decomposition of $\Delta \in \mathcal B(\mathbb T)$, then $\Sigma(\Delta) = \sum_{j=1}^\infty\Sigma(\Delta_j)$ (in the strong sense);
\end{enumerate}
\item An operator-valued measure $\Sigma$ is an  operator-valued probability if $\Sigma(\mathbb T) = I_{\mathfrak H}$;
\item  An operator-valued probability $E$ is said to be spectral or orthogonal if it is projection-valued i.e.  if for any $\Delta \in   \mathcal B(\mathbb T)$, $E(\Delta) = E(\Delta)^2$.
\end{enumerate}
\end{defi}

Let $P_\mathfrak L$ be the orthogonal projection onto the closed linear subspace  $\mathfrak L$. The notation $T_{\!\upharpoonright \mathfrak L}$ means the restriction of a linear operator $T$ on the set $\mathfrak L$.  In the sequel,  $\langle ,  \rangle$ will denote the standard Hermitian product without mention of the ambient Hilbert space. For $\alpha,\beta$ elements of a Hilbert space, the 
outer product $|\alpha\rangle\langle\beta|$ 
 is  a rank one endomorphism defined by 
\[(|\alpha\rangle\langle\beta |)(\cdot)= \alpha \langle\beta, .\rangle\,.\]
%

If $\nu$ is an operator-valued measure on $\mathbb T$ and
  $\mathfrak K$ is a subspace of $\mathfrak H$, then $\nu_{\!\upharpoonright \mathfrak K}$ denotes the map  ${\mathcal B}(\mathbb T) \rightarrow B(\mathfrak K)$ such that $\nu_{\!\upharpoonright \mathfrak K} (\Delta) = \nu(\Delta)_{\!\upharpoonright \mathfrak K}\,.$
Of course, if $\mathfrak L \subset \mathfrak K \subset \mathfrak H$ we have
\ben\label{proj}
\nu_{\!\upharpoonright \mathfrak L } = \left(\nu_{\!\upharpoonright \mathfrak K}\right)_{\!\upharpoonright \mathfrak L}\,.
\een
${\mathcal M}(\mathfrak H)$  (resp. ${\mathcal M}^1(\mathfrak H)$) denotes the set of operator-valued measures in $\mathfrak H$ (resp. operator-valued probability measures in $\mathfrak H$).
We equip  ${\mathcal M}(\mathfrak H)$ with the following topology: $\nu_n \rightarrow \nu$ if, and only if, for all $f \in \mathfrak H$, the sequence of positive (scalar)  measures $\langle f, \nu_n(.)f\rangle$ converges weakly to  $\langle f, \nu(.)f\rangle$.
In the  finite dimensional space case, an operator-valued measure is a matrix measure. For $p\in{\mathbb N}$, we
denote by  ${\mathcal M}_p$  (resp. ${\mathcal M}_p^1$) the set of all Hermitian non-negative $p\times p$ matrix measures  (resp. matrix probability measures).  For $\mu \in \mathcal M_p^1$,  the (matrix) moment $m_\ell(\mu)$ of order $\ell \in \mathbb Z$ of $\mu$ is the element of $\mathbb M_{p,p}$ defined by
\begin{equation}
\label{defmom}m_\ell(\mu) = \int_\mathbb T z^\ell d\mu\,.\end{equation}

If $\lambda$ is the Lebesgue measure on $\mathbb T$, the operator-valued measure $I_p \lambda$ is in some sense the reference measure. We need the notion of absolute continuity and Lebesgue decomposition for operator-valued measures. We refer to Robertson \cite{robertson} and to Mandrekar \cite{mandrekar1}.

  If $\nu$ is a non-negative $\sigma$-finite measure on $\mathbb T$, we say that a $m\times n$ matrix measure $M$ on $\mathbb T$ is absolutely continuous (a.c. in short) with respect to $\nu$ $(M\ll \nu)$ if each entry of $M$ is a.c. with respect to $\nu$. In this case $M'_\nu$ will denote the matrix-valued function whose entries are the Radon-Nikodym derivatives. Let $M$ be a $m\times m$ non-negative Hermitian matrix measure. Then, there exist unique matrix measures $M_a$ and $M_s$ such that $M= M_a + M_s$, $M_a$ and $M_s$ are non-negative. Furthermore,  $M_a \ll \nu$ and $M_s$ is singular with respect to $\nu$, i.e. nonzero only on a set of $\nu$ measure zero.

Now, we say that a matrix measure $Q$ is strongly a.c. with respect to $M$ $(Q\lll M)$ if and only if there exists a nonegative measure $\nu$ dominating $M$ and $Q$ such that the range of the operator $Q'_\nu(\omega)$ is a subset of the range of the operator $M'_\nu(\omega)$  for $\nu$-almost every $\omega$.
Actually, if $M$ is $\ell\times n$ and $Q$ is $m\times n$ then $Q\lll M$ iff there exists a $\ell\times m$ matrix-valued function $\Phi$ integrable with respect to  $M$ such that 
$dQ = \Phi dM$. Such a  $\Phi$ is essentially unique  and  $\Phi = Q'_\nu \left(M'_\nu\right)^\sharp$ where $\sharp$ denotes the pseudo-inverse.

\subsection{Unitary operators 
 and unitary matrices}
\label{suuu}
For any unitary operator $U$ in $\mathfrak H$ the spectral theorem (see \cite{DuSch}) provides  a spectral  operator-valued probability $E_U$  such that
$U = \int_\mathbb T z E_U (dz)$. In other words,  for any $f,g \in \mathfrak H$ and $k \in \mathbb Z$
\ben
\langle f, U^k g \rangle  =  \int_\mathbb T z^k  \langle f, E_U(dz) g\rangle\,.
\een
If $\mathfrak K$ is a subspace of $\mathfrak H$, the spectral measure of the pair $(U, \mathfrak K)$ is by definition  $\left(E_U\right)_{\!\upharpoonright \mathfrak K}$. 
 If $\mathfrak K$ is a one-dimensional subspace, the spectral measure is scalar (see the previous section).

Let $N$ be a fixed integer. In the generic situation,  an operator $U \in \mathbb U(N)$ has  distinct eigenvalues  $ e^{\i\theta_k}, k=1, \cdots, N$ and its  spectral decomposition 
  may be written
\ben
 U = \sum_{k=1}^N e^{\imath\theta_k}\ | v_k\rangle \langle v_k |\,
\,,
 \label{decspc} \een
 where for $k=1,\cdots, N$, $v_k$ is a  unit eigenvector, associated with the eigenvalue $e^{\imath\theta_k}$. 
Obviously,  for $p \leq N$, the spectral measure 
 associated with the pair $(U, \mathcal H_p)$ is
\ben
\label{2.3'}
\mu_{\w, p}:= \sum_{k=1}^N   \w_k  \ \delta_{e^{\imath\theta_k}}\,,\een
 where,  for $k=1,\cdots, N$, $\w_k := \  | P_{\mathcal H_p} v_k \rangle  \langle P_{\mathcal H_p}  v_k |$, is an Hermitian endomorphism on $\mathcal H_p$. 
  This spectral measure is  a matrix probability measure 
which satisfies
\begin{equation}
\label{mom0}
m_\ell(\mu_{\w,p}) = \left(\langle e_i, U^\ell e_j\rangle\right)_{i,j =1}^p\,.
\end{equation}
%
\subsection{Matrix orthogonal polynomials on the unit circle} 
\subsubsection{Construction} 
In spectral theory, orthogonal polynomials play a prominent role. Here, we will need to work with matrix-valued orthogonal polynomials with respect to matrix measures on the unit circle. We recall some useful facts and refer to \cite{simon05} for more on the subject.
To begin with, as in the scalar case, these orthogonal polynomials satisfy a recursion and the matrices appearing in this recursion are the so-called matrix Verblunsky coefficients (see \cite{damanik2008analytic}). Let us give some more notations. Let $p\in\mathbb{N}$  and $\mu\in \mathcal M_p$. 
Further, let  $F$ and $G$ be 
measurable matrix valued  functions  : $\mathbb T \rightarrow \mathbb M_{p,p}$. We define two $p\times p$ matrices by setting
\begin{eqnarray*}\lal F, G\rar_R &=& \int_\mathbb T F(z)^\dagger d\mu(z)G(z) \in \mathbb M_{p,p}\\
\lal F, G\rar_L &=& \int_\mathbb T G(z) d\mu(z)F(z)^\dagger \in \mathbb M_{p,p}\,.
\end{eqnarray*}
A sequence of functions  $(\varphi_j)$ in $\mathcal H_p$  is called right-orthonormal if and only if 
\[\lal\varphi_i, \varphi_j\rar_R= \delta_{ij}I_p \ .\]
The orthogonal polynomial recursion is built as follows. 

First, assume that the support of $\mu$ is infinite. We define the right monic matrix orthogonal polynomials $\Phi_n^R$ by applying Gram-Schmidt procedure  to $\{I_p, zI_p, z^2 I_p,\dots\}$.  In other words, $\Phi_n^R$ is the unique matrix polynomial $ \Phi_n^R(z) = z^nI_p +$ lower order terms such that $\lal z^kI_p, \Phi_n^R\rar_R =0$ for $k=0, \dots, n-1$. 
The normalized orthogonal polynomials are defined by
\[\varphi_0 = I_p\ \  ,\ \  \varphi_n^R = \Phi_n^R\kappa_n^R\]
where the sequence of $p\times p$ matrix $(\kappa_n^R)$ 
satisfy, for all $n$, the condition $\left(\kappa_n^R\right)^{-1}\kappa_{n+1}^R > 0$
and is such that the set $\{\varphi_n^R\}$ is orthonormal.
We define the sequence of left-orthonormal polynomials $\{\varphi_n^L\}$ in the same way except that the above condition is replaced by  $\kappa_{n+1}^L \left(\kappa_n^L\right)^{-1}> 0$. The Szeg\H{o} recursion is then 
\begin{eqnarray}
\label{SzL}
z\varphi_n^L -\rho_n^L\varphi_{n+1}^L &=& \alpha_n^\dagger(\varphi_n^R)^*\\
\label{SzR}
z\varphi_n^R - \varphi_{n+1}^R\rho_n^R &=& (\varphi_n^L)^* \alpha_n^\dagger\,,
\end{eqnarray}
where for all $n\in\mathbb{N}$,  
\begin{itemize}
\item $\alpha_n$ belongs to  $\mathbb{B}_p$ the closed unit ball of $\mathbb M_{p\times p}$  defined by
$$\mathbb{B}_p:=\{V \in \mathbb M_{p,p} : VV^\dagger \leq I_p\}\,,$$
\item $\rho_n$ is the so-called defect matrix defined by
\ben
\label{defrho}
\rho^R_n :=  \left(I_p - \alpha_n\alpha_n^\dagger\right)^{1/2}\  , \  \rho^L_n =    \left(I_p - \alpha_n^\dagger\alpha_n\right)^{1/2}\,,
\een
\item for a matrix polynomial $P$, having degree $n$,  the reversed polynomial $P^*$ is defined by
\[P^*(z) := z^n P(1/\bar z)^\dagger\,.\]
\end{itemize}
Notice that  the construction of the recursion coefficients uses only the matrix moments.
The Verblunski's theorem (analogue of Favard's theorem for matrix orthogonal polynomials on the unit circle) establishes a one-to-one correspondance between matrix measures on $\mathbb T$ with infinite support and sequences of elements of the interior of $\mathbb B_{p}$ (Theorem 3.12 in \cite{damanik2008analytic}). 

Now, for a matrix measure having a finite support, the construction of the Verblunsky coefficients is not obvious. In \cite{dewag10} Theorem 2.1, a sufficient condition on the moments for such a construction is provided. It is related to the positivity of a block-Toeplitz matrix, as it is also mentioned  in \cite{simon05} at the top of 
  p.208.

\begin{lem}
\label{lemdur}
Let $N = Qp+r$, with $Q \geq 1$ and $0 \leq r<p$. Let $J\leq Q-1$. 
 Then for almost every $U\in \mathbb U(N)$ equipped with the Haar distribution,
 there exists a measure $\nu$ with infinite support and satisfying
\[m_k(\nu) = m_k(\mu_{\w,p}) \ , \  k \leq J\,.\]
\end{lem}

\subsubsection{A unitary isomorphism}
As in the scalar case, we can build a unitary equivalence between the linear space $\mathbb M_{N,p}$ (of systems of $p$ vectors of ${\mathcal H}_N$) and the linear space $\mathbb P_{N,p}$ of polynomials of degree $\leq N-1$ with coefficients in $\mathbb M_{p,p}$. 
Let $e = [e_1 \cdots e_p]$ be  a $N \times p$ matrix consisting in $p$ column vectors of dimension $N$. In the same way, we consider the $N \times p$ matrix $Ue = [Ue_1 \cdots Ue_p]$. The pseudo-scalar (inner) product of $e$ with  $f = [f_1 \cdots f_p]$ is a matrix $p\times p$ denoted by $\ll f, e\gg$ and  defined by
\[\ll f, e\gg_{i,j}\!\ = \langle f_i , e _j\rangle \ \  i,j = 1, \dots, p\,.\]
It is clear that if the system $e_1, \cdots ,e_p$ is orthonormal, then $\ll e, e\gg = I_p$. It is clear also that if $W\in \mathbb M_{N, N}$
 and $e,f$ are as above, then
\begin{equation}
\label{adj}
\ll f, We\gg \!\ =\!\ \ll W^\dagger f, e\gg\,.\end{equation} 
Elementary computations lead also to the useful following properties
\begin{equation}
\label{magicllrr}
\ll f,gW\gg \!\ =\!\ \ll f,g\gg W \  , \ \ll fW,g\gg\!\ =\!\ W^\dagger\ll f,g\gg \,,
\end{equation}
where $W \in \mathbb M_{p,p}$.
Notice that  $s = [s_1 \cdots s_p]$  consists in  elements of  $\spann \{e\}$ if and only if there exists a matrix $\gamma\in \mathbb M_{p,p}$ such that $s = e\gamma$.

Let $\varepsilon =[\varepsilon_1\cdots\varepsilon_p]$ an orthonormal system, $U$ unitary on $\mathcal H_N$ and $\mu$ the spectral measure of the pair $(U, \spann\varepsilon)$.
\begin{defi}
We say that  $\varepsilon$ is cyclic for $U$, if 
\ben
\spann \{ U^k \varepsilon_j , 0 \leq k \leq N-1, 1\leq j\leq p\} = \mathbb C^N\,. 
\een
\end{defi}
In this case, each element $M$ of $\mathbb M_{N,p}$ may be written as 
$M= \sum_{k=0}^{N-1} U^k\varepsilon\gamma_k$ where $\gamma_k \in \mathbb M_{p,p}$  and then we can associate the polynomial $q(z) = \sum_{k=0}^{N-1} \gamma_k z^k \in \mathbb P_{N,p}$. 

\begin{rem}
The system  $\varepsilon = [\varepsilon_1, \dots, \varepsilon_p]$ is cyclic for $U$ 
as soon as one of the $\varepsilon_j$ is cyclic for $U$, but of course, it is not a necessary condition.
\end{rem}
 If $\varepsilon$ is cyclic for $U$ we have then a one-to-one correspondance  between $\mathbb M_{N,p}$ and $\mathbb P_{N,p}$ which preserves the pseudo-scalar products $\ll , \gg$ and $\langle\langle , \rangle\rangle_R$:
\begin{itemize}
\item To $q(z) = \sum_k \gamma_k z^k \in \mathbb P_{N,p}$ we associate 
 $M_q = \sum_k U^k\varepsilon \gamma_k \in \mathbb M_{N,p}$
\item If $p$ and $q$ are two polynomials, we have
\[ \ll M_p, M_q\gg\!\ = \langle\langle p, q\rangle\rangle_R\,,\]
(this property  is straightforward for monomials and is extended easily). 
\end{itemize} 
\begin{lem}
\label{lemfou}
If $\varepsilon$ is a cyclic system for $U$, then the first Verblunsky coefficient of the pair $(U, \varepsilon)$ denoted by $\alpha_0(U, \varepsilon)$  satisfies
\begin{equation}
\label{relalpha0U}
\alpha_0(U, \varepsilon)^\dagger = \!\ \ll \varepsilon, U\varepsilon\gg\,. \end{equation}
\end{lem}
We construct now the GGT matrix $\mathcal G^R$, which is the matrix of the unitary operator 
$$f \in L^2(d\mu) \mapsto \left(z \mapsto zf(z)\right) \in L^2(d\mu)$$ in the orthonormal basis $\varphi_k^R$, that is
\[\mathcal G^R_{k, \ell} = \langle\langle\varphi_k^R , z \varphi^R_\ell\rangle\rangle_R\,. \]

\begin{lem}
\label{lemCGT}
We have
\begin{eqnarray}
{\mathcal G}_{k,\ell}^R = \begin{cases} - \alpha_{k-1}\rho_k^L\rho_{k+1}^L\cdots \rho_{\ell -1}^L \alpha_\ell^\dagger \;\;\;\;& 0\leq k\leq\ell\\
\rho_\ell^R &k= \ell+1\\
0&k \geq \ell+1
\end{cases}
\end{eqnarray}
(with $\alpha_{-1} = -1_p\ , \ \rho^L_{-1}=0$).
\end{lem}
It is convenient to call $\mathcal G^R(\alpha_0, \alpha_1, \dots)$ the above GGT matrix built from the Verblunsky coefficients $(\alpha_0, \alpha_1, \dots)$. We have the easy following result, which is a replica of  Theorem 10.1 in  Simon  \cite{simon2006cmf}.
\begin{prop}
Let $\alpha$ be in the matrix unit ball. Set \ben
\label{deftheta0}
\Theta(\alpha):=
\begin{pmatrix}
\alpha^{\dagger}&\rho^L\\
\rho^R& -\alpha
\end{pmatrix}\,,
\een 
and $\widetilde \Theta(\alpha) := \Theta(\alpha) \oplus I_p\oplus I_p \oplus...$.
 Then we have
\begin{equation}
\label{AGRGGT}\mathcal G^R\left(\alpha_0, \alpha_1, \dots\right) = \widetilde \Theta(\alpha_0)[I_p \oplus \mathcal G^R\left(\alpha_1, \alpha_2, \dots\right)]\,.\end{equation}
\end{prop}

For the sake of completeness, let us consider now the operator point of view.
\begin{prop}
\label{new1}
Assume that $\varepsilon$ is cyclic for $U$. 
Let $\xi(U, \varepsilon) \in \mathbb M_{N,p}$ be the image of $\varphi_1$ in the isomorphism and set $H_0 := \spann\{\varepsilon, U\varepsilon\} = \spann\{\varepsilon, \xi\}$. Let $V(U, \varepsilon)$ be the unitary transform which leaves invariant the subspace orthogonal to 
$H_0$ and whose restriction to $H_0$ has the matrix $\Theta(\alpha_0(U, \varepsilon))$ in the basis $(\varepsilon, \xi)$.
Then $\xi(U, \varepsilon)$ is cyclic for the restriction of $W(U, \varepsilon) :=  V^{-1}(U, \varepsilon)U$ to $\varepsilon^\perp$ and we have
\ben\label{shift}
\alpha_1(\varepsilon, U) = \alpha_0(\xi(U, \varepsilon), W(U, \varepsilon))\,.
\een
\end{prop}

\section{Randomization}
\label{srando}
\subsection{Distributions on ${\mathbb M}_{n,n}$}
\label{suDistri}
Along this paper, we use three probability distributions,
\begin{itemize}
\item The Haar measure on $\mathbb U (n)$ for $n \geq 1$.
\item $\tt Gin$$(n)$, the (Ginibre) distribution on  ${\mathbb M}_{n,n}$ which makes all elements independent and standard complex Gaussian. It has the density
\[\pi^{-n^2} \exp \left(- \tr\!\ GG^\dagger\right)\]
with respect to the Lebesgue measure.
\item $\tt Cor$$(n,p)$ (for$\  n > 2p$) the distribution on  ${\mathbb M}_{p,p}$ of the top-left corner of size $p\times p$ of a Haar distributed random matrix. It has the density
\ben
\label{33}V \mapsto K_{p,n} \det\left(I_p - V V^\dagger\right)^{n-2p}\een
on the unit ball $\mathbb{B}_p$ where $K_{p,n}$ is the normalization constant (see for instance \cite{collins2005product} Theorem 5.1). Actually,  \cite{Ner1} Lemma 1.4 gives
\[K_{p,n} = \pi^{-p^2}\frac{(n-2p)! (n-2p+1)! \dots (n-p-1)!}{(n-p)!(n-p+1)!\dots (n-1)!}.\]
\end{itemize}
\subsection{Preliminary results}
\label{suprere}
Here, we describe the distribution of the matrix spectral measure.  
The first statement uses the encoding by weights and (\ref{2.3'}).
\begin{prop}
\label{propsursoi}
Let $U$ be   Haar distributed in $\mathbb U(N)$.
\begin{enumerate}
\item
The random variable $(e^{\i \theta_1}, \cdots, e^{\i \theta_N})\in \mathbb T^N$ is independent of the random variable $\left(\w_1, \cdots, \w_N\right)$.
\item  For $ k=1, \cdots, N$,  let $a_k$  be independent $p$-dimensional random vectors with complex standard normal distribution.  
The random variable 
 $\left(\w_1, \cdots, \w_N\right)$ has the same distribution as
\[\left(\h^{-1/2}\v_1 \h^{-1/2}, \cdots, \h^{-1/2}\v_N \h^{-1/2}\right)\]
where $\v_k := | a_k\rangle a_k |\ , \ k=1, \cdots, N$ and 
\ben
\label{H=} \h := \sum_{k=1}^N \v_k .\een

\end{enumerate}
\end{prop}
The second statement  describes the distribution of the  Verblunsky coefficients.

\begin{theo}
\label{alphajlaw}
Let $N = pQ +r $ with $0 \leq  r < p$ and $Q > 2$. Let $U$ be chosen at random in $\mathbb U(N)$ according to Haar measure. Let $\mu_p\sn$ be the matrix spectral measure for 
$(U, e_1, \dots, e_p)$. 
 Then, the  matrix Verblunsky coefficients 
$\alpha_j\sn := \alpha_j(\mu_p\sn)$ for $ j=0, \cdots, Q-2$
are independent. Moreover, for $j \leq Q-2$, 
$\alpha_j\sn$ has in the matrix unit ball of $\mathbb M_{p,p}$ the density $\tt Cor$$(N-pj, p)$.
%
\end{theo}

\begin{rem}
Notice that for $p=1$ the last Verblunsky coefficient $\alpha_{N-1}^{(1)}$ is uniformly distributed on $\mathbb T$. In the general case, the distribution of what could be the last coefficient $\alpha_{Q-1}\sn$ is not obvious.  
Assume that  $N= p+r$ with $r < p$. Let us compute $\alpha_1$.  Let us denote
\[U = \begin{pmatrix} \alpha_0 & C\\B&A
\end{pmatrix}\,.\]
According to Arlinskii \cite{Arlinski} Theorem 4.2, we have
\begin{equation}\label{arl}\alpha_1 = (\rho_0^R)^{-1}CB(\rho_0^L)^{-1}\end{equation}
where $(\rho_0^R)^{-1}$ and $ (\rho_0^L)^{-1}$ are the  Moore-Penrose pseudo-inverses of  $\rho_0^R$ and $ \rho_0^L$ respectively. To see that we cannot use the true inverses let us look at
 $(\rho_0^R)^2 = I_p - \alpha_0\alpha_0^\dagger= CC^\dagger$. 
The singularity of $\rho_0^R$ follows directly from
\[\hbox{rank}\ (CC^\dagger) \leq \hbox{rank}\ (C) < p\,. \]
 We did not succeed in computing the distribution of $\alpha_1$ given by (\ref{arl}) in such a case. Notice that in formula (\ref{arl}) we can replace $CB$ by $\Gamma_2 - \Gamma_1^2$ where $\Gamma_1 = \alpha_0$ and $\Gamma_2$ are the first two moments of $\mu$. We recover a formula which fits with (2.19) in \cite{dewag10}.

\end{rem}
\subsubsection{Asymptotics}
In this section, we consider central limit theorems (CLT) i.e. convergences in distribution, and large deviation principles (LDP). To make this paper self-contained, let us first recall what is a LDP. For more on LDP we refer to \cite{demboZ98}.

 Let $(a_N)$
be an increasing positive sequence of real numbers going
to infinity with $N$.
\begin{def}\label{dldp}
We say that a sequence $(Q_N)$ of probability measures on a measurable
Hausdorff space $U$ with corresponding Borel field $\br(U))$ {\rm satisfies an LDP at scale $a_N$
and rate function $I$ } if:
\begin{itemize}
\item[i)] $I$ is lower semicontinuous, with values in
$\R^{+}\cup\{+\infty\}$.
\item[ii)] For any measurable set $A$ of $U$:
$$-I(\Int A)\leq
\liminf_{N\rightarrow\infty}\usn^{-1}\log Q_N(A)\leq
\limsup_{N\rightarrow\infty}\usn^{-1}\log Q_N(A)\leq
-I({\Clo A}),$$
where $I(A)=\inf_{\xi\in A}I(\xi)$ and $\Int A$ (resp. $\Clo A$) is the
interior (resp. the closure) of $A$.
\end{itemize}%
We say that the rate function $I$ is {\rm good} if its level sets
$\{x\in U:\; I(x)\leq a\}$ are compact for any $a\geq 0$.
More generally, a sequence of $U$-valued random variables
is said to satisfy an
LDP if their distributions satisfy a LDP.
\end{def}
\begin{prop} [CLT] 
$\;$
\label{CLT}
\begin{enumerate}
\item If  $V_p\sn$ is a random matrix with distribution $\tt Cor$ $(N,p)$, 
 then for fixed $p$
\[{\sqrt N} V_p\sn \schw {\tt Gin} (p)\]
\item
 If $k$ and $p$ are fixed, the $k$ first (matrix) Verblunsky coefficients satisfy
\[{\sqrt N}(\alpha_0\sn, \dots, \alpha_k\sn)  \schw (G_1, \dots, G_k)\]
where $G_1, \dots ,G_k$ are independent and  $\tt Gin$$(p)$ distributed.
\end{enumerate}
\end{prop}
There are two proofs of the first statement in \cite{118218FGAR}. The second statement is a consequence of the first statement  and of Theorem \ref{alphajlaw} above.

Coming back to the moments of the measure $\mu_\w\sn$  we recover a result of Krishnapur,  
so offering a proof (postponed to Section \ref{proofharry}) shorter than the involved combinatorial original one.

\begin{cor}[\cite{krishna09},  Lemma 10 p.357]
\label{corharry}
Let $U \in \mathbb U(N)$  sampled from the Haar measure. Fix $p \geq 1$ and $n_0 \geq 1$.  Then the sequence of random variables  $\sqrt N [U^n]_{i,j},  1 \leq n \leq n_0, i,j \leq p$, converges as $N \rightarrow \infty$ in distribution to independent standard complex Gaussian matrices. 
\end{cor}

\begin{theo}[LDP]
\label{ldpalpha}
$\;$
\begin{enumerate}
\item
 If  $V_p\sn$ is a random matrix with density (\ref{33}) in the unit ball, then for fixed $p$, the sequence $(V_p\sn)_N$ satisfies the LDP at scale $N$  in $\mathbb M_{p,p}$ with good rate function
\begin{align}
v \in \mathbb M_{p,p} \longmapsto I(v) =
\begin{cases}
 -\log\det (I_p- vv^\dagger)\;\;\;\;&\mbox{if}\;vv^\dagger < I_p,\\
      \infty\;\;&\mbox{ otherwise.}
\end{cases}
\end{align} 
\item
For fixed $p\geq 1$ and $k\geq 0$, the sequence $\left(\alpha_0\sn, \dots, \alpha_k\sn\right)_N$ satisfies the LDP at scale $N$ in $\mathbb M_{p,p}\times \cdots \times\mathbb M_{p,p}$ with good rate function 
\begin{align}I^k (\alpha_0, \dots, \alpha_k) = \begin{cases}- \sum_{j=0}^k \log \det (I-\alpha_j^\dagger\alpha_j)\;\;\;\;&\mbox{if}\;\alpha_j^\dagger\alpha_j < I_p\; \mbox{for}\;  j=0, \dots, k\,,\\
      \infty\;\;&\mbox{ otherwise.}
\end{cases}
\end{align} 
\end{enumerate}
\end{theo}

The first statement is a direct consequence of the explicit expression of the density and the second statement comes from the independence of Verblunsky coefficients. These are arguments from \cite{118218FGAR}. It is worthwhile to quote the scalar case which was established in \cite{LozEJP}.

\subsection{Large deviations for the spectral measure in fixed dimension}
\begin{theo}
\label{rateptheo}
For fixed $p \geq 1$, the family of random matrix measures $\left(\mu_{\w, p}\sn\right)_N$ satisfies the LDP in ${\mathcal M}_p^1$ at scale $N$ 
with good rate function
\ben
\label{ratep}
\nu \in \mathcal M_p^1 \longmapsto {\mathcal I}_p (\nu) = -  \int_\mathbb T \log \det\nu'_a(z) dz
\een
where $d\nu(z) = \nu'_a(z) dz + d\nu_s(z)$ is the Lebesgue decomposition of $\nu$.
\end{theo}
\begin{cor}
\label{rem}
It is possible to rewrite the above quantity in the flavour of Kullback information with the notation of \cite{mandrekar1} or \cite{robertson}, i.e.
\ben
\label{ratep0}
{\mathcal I}_p (\nu) = 
\begin{cases}
\displaystyle \int_\mathbb T \log \det \frac{I_p dz}{d\nu(z)} \!\ dz \;\;\;\;&\mbox{if}\; I_p dz \lll d\nu(z),\\
      \infty\;\;&\mbox{ otherwise.}
\end{cases}
\een
where  $\lll$ denotes the strongly absolute continuity (see Section \ref{suopsp}).
\end{cor}

\begin{rem}
\label{rem1}
The deviations are from $\lambda_p(dz) :=I_p dz$ whose  moments of every order  are zero, i.e. $\int_\mathbb T z^k I_p dz = 0_p$ for $k \not=0$, where $0_p$ is the null endomorphism on $\mathcal H_p$. 
This corresponds to the fact that $\lim_N (U^k)_{i,j} =0$ for every $k,i,j\geq 1$ fixed.
%
\end{rem}
\subsection{Large deviations - Operator-valued random measures}

Every element $U$  of $\mathbb U(N)$ can be extended to an operator on $\mathcal H$ by tensorisation with identity.

More precisely, if $(e_i)$ denotes as above the canonical basis of ${\mathcal H}$ we define
\[\widehat U\left(\sum_1^\infty h_i e_i\right)  = \sum_1^N h_i U(e_i) + \sum_{N+1}^\infty h_j e_j\,.\]
When $U$ is chosen according to the Haar measure in $\mathbb U (N)$, the (random) spectral measure associated with  $\widehat U$ is denoted $\mu\sn$. It is of the form
\[\mu\sn = E_{\widehat U} = \sum_{k=1}^N |\widehat v_k\rangle \langle \widehat v_k | \!\  \delta_{e^{i\theta_k}} + \left(\sum_{k=N+1}^\infty | e_k\rangle\langle e_k | \right) \delta_1\,,\]
where the $\widehat v_k$ are the eigenvectors of $U$ extended in $\mathbb C^N$ by zeros.
We establish now an LDP for this sequence.

\begin{theo}
\label{theoarem}
The family of random 
  spectral
  measures $\left(\mu\sn\right)_N$ satisfies the LDP in ${\mathcal M}^1(\mathcal H)$ at scale $N$ 
with the good rate function  
\begin{equation}
\label{iinfini}
\mu \in {\mathcal M}^1(\mathcal H)\longmapsto {\mathcal I}_\infty (\mu) = \lim_k\uparrow\int_\mathbb T - \log \det {\mu'_a}^k (z)
 dz\end{equation}
where  
$\mu = \mu'_a dz + \mu_s$
 and $\mu'_a$  
 is a measurable function from $\mathbb T$ to 
$H(\mathcal H)$ and 
\[{\mu'_a}^k  
 = P_{{\mathcal H}_k} \left(
\mu'_a(z)
\right)_{\!\upharpoonright {\mathcal H}_k}\,.\] Moreover, if there is a constant $C > 1$ such that for every $k$ and $z$
\begin{equation}
\label{C}C^{-1} \leq \det
{\mu'_a}^k(z) \leq C\end{equation}
and if
 for every $z$ the operator  $I_\mathcal H -\mu'_a(z)$ 
is trace class and $z \mapsto \tr \log 
\mu'_a(z) \in L^1(\mathbb T)$, then 
\ben
\label{ratepinfty}
{\mathcal I}_\infty (\mu) = 
\displaystyle \int_\mathbb T - \tr \log  
\mu'_a(z) dz\,. 
\een
\end{theo}
\begin{rem}
\label{remrem1}
The deviations are from $\lambda_\infty(dz) :=I_\mathcal H dz$, whose  moments of every order  are zero, i.e.  $\int_\mathbb T z^k I_\mathcal H dz =0_\mathcal H$
 for $k \not=0$, where $0_\mathcal H$ is the null endomorphism on $\mathcal H$. As in Remark \ref{rem1} this corresponds to the fact that $\lim_N (U^k)_{i,j} =0$ for every $k,i,j\geq 1$ fixed. 
\end{rem}
\begin{rem}
\label{remrem2}
The boundedness assumption on the sequence $(\det
{\mu'_a}^k(z))$ cannot be deduced from a simple hypothesis on the operator valued density $\mu'_a$. 
\end{rem}

\section{Proofs}
\label{sproofs}
\subsection{Proof of Lemma \ref{lemdur}}
Let us first describe  the result of Dette and Wagener \cite{dewag10} (up to  a slight change of notation). Let for $\ell = 1,2, \dots$
\[F_\ell := \{(I_p, m_1(\mu), \dots, m_\ell(\mu))\  | \ \mu \in {\mathcal M}_p^1\}\]
the $\ell$th moment space of matrix (probability) measures on the unit circle.
Their Theorem 2.1 says that 
$\widetilde m = (1, m_1, \dots, m_J) \in F_J$ if and only if 
\[\sum_{i,j =0}^J \tr (B_j^\dagger m_{i-j}B_i) \geq 0 \ \ \hbox{for  all} \ \ B_0, \dots, B_J \in \mathbb C ^{p\times p}\]
where for $k > 0, \  m_{-k} = m_k^\dagger$. Moreover $\widetilde m
 \in \hbox{int}\  F_J$ if and only if there is a strict inequality above except if $B_0= \dots= B_\ell = 0$.

For $\ell =1, \dots, J$, let  $m_\ell := m_\ell (\mu_{\w, p})$. We have 
\begin{eqnarray}\sum_{i,j =0}^J \tr (B_j^\dagger m_{i-j}B_i) &=& \sum_{k=1}^N \tr \left(\sum_{i,j =0}^J B_j^\dagger P_{\mathcal H_p} v_k (P_{\mathcal H_p} v_k)^\dagger B_i e^{\imath\theta_k(i-j)}\right)\\
\nonumber
&=&  \sum_{k=1}^N \tr \ A_k^\dagger A_k\,,
\end{eqnarray}
where $A_k = \sum_{j=0}^J e^{\imath j \theta_k} (P_{\mathcal H_p} v_k)^\dagger B_j$.
Let us prove that in our conditions $(1, m_1, \dots, m_J) \in \hbox{int}\ F_J$ a.s. That means that a.s. we cannot find a system $(B_0, \dots, B_J)$ non zero such that 
\begin{equation}\label{sumak}
\sum_{k=1}^N \tr \ A_k^\dagger A_k =0\,.
\end{equation} But since the matrices $A_k^\dagger A_k$ are Hermitian nonnegative, (\ref{sumak}) is equivalent to $A_k = 0$ for all $k= 1, \dots, N$ which yields to a system of $N$ matricial equations:
\[ \sum_{j=0}^J e^{\imath j \theta_k} (P_{\mathcal H_p} v_k)^\dagger B_j = 0 \ , \ \  k= 1, \dots, N\,,\]
which may be converted into  a linear system of $pN$ scalar equations in the $p^2 (J+1)$ unknown variables $(B_j)_{s,t}$ where 
$j = 0, \dots ,J$ and  $s,t =1, \dots, p$.
We have supposed that $J \leq Q-1$, so that in all cases $p^2(J+1) \leq pN$. 
If this system has a nontrivial solution, 
the following $p^2(J+1)\times p^2(J+1)$ minor 
\[\begin{pmatrix}
I_p\otimes \Gamma_1\\
\dots\\
I_p \otimes \Gamma_{(J+1)p}
\end{pmatrix} \ \hbox{where}\ \Gamma_k = \left( P_{\mathcal H_p} v_k^\dagger ,  P_{\mathcal H_p} v_k^\dagger e^{\imath\theta_k},  P_{\mathcal H_p} v_k^\dagger e^{2\imath\theta_k},  \dots,  P_{\mathcal H_p} v_k^\dagger e^{J\imath\theta_k}\right)\]
has a determinant zero. An easy permutation of rows shows that this determinant is, up to a change of sign,
\[\left[\det \begin{pmatrix}
\Gamma_1\\
\vdots\\
\Gamma_{(J+1)p}
\end{pmatrix}\right ]^p
\]
If our system has a non-trivial solution, we have
\begin{equation}
\label{algebrique}
\det \begin{pmatrix}
\Gamma_1\\
\vdots\\
\Gamma_{(J+1)p}
\end{pmatrix} =0\,.
\end{equation}

This equation is polynomial in the variables $e^{\imath\theta_\ell},\left(P_{\mathcal H_p}v_{k}\right)^T$ for $\ell, k \leq (J+1)p$. 
Besides, under the Haar probability on $\mathbb U (N)$ the distribution of $e^{\imath\theta_\ell}, \ell \leq (J+1)p$ is absolutely continuous with respect to the Lebesgue measure on $\mathbb T^{Jp}$,  the distribution of $\left(P_{\mathcal H_p}v_{k}\right)^T , k \leq (J+1)p$ is also diffuse on the set $\{\sum_{s=1}^{(J+1)p} f_s^Tf_s \leq I_p\}$ and both arrays are independent.  If the above  polynomial is not identically zero, the set of its solutions will then be of measure zero. But if we choose 
$P_{\mathcal H_p}v_1, \dots , P_{\mathcal H_p}v_p$ as the canonical basis $(e_1, \dots, e_p)$ of $\mathbb C^p$ and $P_{\mathcal H_p}v_{kp+ s} = P_{\mathcal H_p}v_{s}= e_s$ for $k= 1, \dots, Jp$ and $s =1, \dots, p$, equation (\ref{algebrique}) becomes:
\begin{equation}
\label{VDM}
\det \begin{pmatrix}
e_1^T&e^{\imath\theta_1}e_1^T& e^{2\imath\theta_1}e_1^T&  \dots& e^{J\imath\theta_1}e_1^T\\
e_2^T&e^{\imath\theta_2}e_2^T& e^{2\imath\theta_2}e_2^T&  \dots & e^{J\imath\theta_2}e_2^T\\
\vdots&\vdots&\vdots&\vdots&\vdots\\
e_p^T&e^{\imath\theta_p}e_p^T& e^{2\imath\theta_p}e_p^T& \dots & e^{J\imath\theta_p}e_p^T\\
e_1^T&e^{\imath\theta_{p+1}}e_1^T& e^{2\imath\theta_{p+1}}e_1^T& \dots & e^{J\imath\theta_{p+1}}e_1^T\\
\vdots&\vdots&\vdots&\vdots&\vdots\\
e_p^T& e^{\imath\theta_{2p}}e_p^T& e^{2\imath\theta_{2p}}e_p^T& \dots & e^{J\imath\theta_{2p}}e_p^T\\
\vdots&\vdots&\vdots&\vdots&\vdots\\
e_1^T&e^{\imath\theta_{Jp+1}}e_1^T& e^{2\imath\theta_{Jp+1}}e_1^T&\dots & e^{J\imath\theta_{Jp+1}}e_1^T\\
\vdots&\vdots&\vdots&\vdots&\vdots\\
e_p^T&e^{\imath\theta_{(J+1)p}}e_p^T& e^{2\imath\theta_{(J+1)p}}e_p^T&\dots & e^{J\imath\theta_{(J+1)p}}e_p^T\\
\end{pmatrix} =0\,.
\end{equation}

Afer reordering into Vandermonde-type blocks,
 the right hand side is (up to a sign)
$$\prod_{k=1}^p\prod_{0\leq i <j\leq J} \left(e^{\imath\theta_{ip+k}} - e^{\imath\theta_{jp+k}}\right)\,.$$
So, it is not identically zero and we may conclude that a.s. $\left(I_p, m_1(\mu_{\w,p}), \dots , m_J(\mu_{\w,p})\right) \in \ \hbox{int}\ F_p$. Hence, we can construct a measure $\nu$ with infinite support whose $J$ first moments fit with those of $\mu_{\w, p}$ using the Bernstein-Szeg\H{o} construction (see \cite{damanik2008analytic} Section 3.6).

\QED
\subsection{Proof of Lemma \ref{lemfou}}
We start with $\varphi_0^R = \varphi_0^L= I_p$. 
Then by (\ref{SzL}) and  (\ref{SzR}) with $n=0$
\begin{eqnarray}
\nonumber
zI_p  - \rho_0^L\varphi_1^L &=& \alpha_0^\dagger\\
\label{Sz0}
zI_p - \varphi_1^R\rho_0^R &=& \alpha_0^\dagger\,.
\end{eqnarray}
Since $\varphi_1^R$ (resp. $\varphi_1^L$)  is orthogonal to $\varphi_0^R$ (resp. $\varphi_0^L$) we get
$\alpha_0^\dagger = \langle\langle I_p, zI_p\rangle\rangle_R$,  
so that
$\alpha_0(U, \varepsilon)^\dagger = \!\ \ll \varepsilon, U\varepsilon\gg$.
\QED
\subsection{Proof of Lemma \ref{lemCGT}}
In this proof, all inner products are right inner products.
Let us begin with the subdiagonal terms:
\[\mathcal G_{\ell +1, \ell} = \langle\langle \varphi_{\ell +1}^R, z\varphi_\ell^R\rangle\rangle\,.\]
From (\ref{SzR}))
\[\langle\langle \varphi_{\ell +1}^R, z\varphi_\ell^R\rangle\rangle = \langle\langle \varphi_{\ell +1}^R, \varphi_{\ell +1}^R\rangle\rangle \rho_\ell^R
- \langle\langle \varphi_{\ell +1}^R, (\varphi_\ell^L)^*\rangle\rangle \alpha_\ell^\dagger\]
and since $ (\varphi_\ell^L)^*$ is a polynomial of degree $\ell$, it is orthogonal to $ \varphi_{\ell +1}^R$, and 
$\mathcal G_{\ell+1, \ell} = \rho_\ell^R$.

Suppose now $0 \leq k \leq \ell$. 
Again from (2.6) we have
\[ \langle\langle\varphi_k^R , z \varphi^R_\ell\rangle\rangle = \langle\langle \varphi_k^R ,  \varphi_{\ell +1}^R\rangle\rangle\rho_\ell^R -  \langle\langle \varphi_k^R   , (\varphi_\ell^L)^*\rangle\rangle \alpha_\ell^\dagger\]
Now it is the first term which vanishes, so that it remains 
\ben
\label{firstterm}
 \langle\langle\varphi_k^R , z \varphi^R_\ell\rangle\rangle =  -  \langle\langle \varphi_k^R   , (\varphi_\ell^L)^*\rangle\rangle \alpha_\ell^\dagger
\een
In the Christoffel-Darboux formula (see \cite{damanik2008analytic} Proposition 3.6 (b))
\[(1-\bar z_1 z_2) \sum_{k=0}^n \varphi_k^R(z_2) \varphi_k^R(z_1)^\dagger = \varphi_{n+1}^{L, *}(z_2) \varphi_{n+1}^{L,*}(z_1)^\dagger - \varphi_{n+1}^{R}(z_2) \varphi_{n+1}^R(z_1)^\dagger \]
taking $z_1=0, z_2=z$ gives
\ben
\varphi_\ell^{L,*}(z) = \sum_{k=0}^\ell \varphi_k^R(z)g_k^\ell \ , \ g_k^\ell = \varphi_k^R(0)^\dagger\left(\varphi_\ell^{L,*}(0)^\dagger\right)^{-1}
\een
and then
\ben
\label{phiphistar}
 \langle\langle \varphi_k^R   , (\varphi_\ell^L)^*\rangle\rangle = g_k^\ell
\een
But we have 
\ben
\varphi_k^R = \Phi_k^R\kappa_k^R\ ; \ \varphi_\ell^L = \kappa_\ell^L\Phi_\ell^L \ ; \  \Phi_\ell^{L,*}(0) = 1_p
\een
and from the Szeg\H{o} recursion
\[\varphi_k^R(0) = - \left(\kappa_{k-1}^L\right)^\dagger \alpha_{k-1}^\dagger\left(\rho_{k-1}^R\right)^{-1}\]
so that
\ben
\nonumber
g_k^\ell = - (\rho_{k-1}^R)^{-1} \alpha_{k-1}\kappa_{k-1}^L (\kappa_\ell^L)^{-1}&=& -  \alpha_{k-1} (\rho_{k-1}^L)^{-1}\kappa_{k-1}^L (\kappa_\ell^L)^{-1}\\
\label{gkl}&=& -\alpha_{k-1} \kappa_k^L (\kappa_\ell^L)^{-1}
\een
(we have used  $\alpha_j\rho_j^L = \rho_j^R\alpha_j$ and  $\rho^L_{j-1} = \kappa_{j-1}^L\left(\kappa_j^L\right)^{-1}$).
Gathering (\ref{gkl}), (\ref{phiphistar}) and (\ref{firstterm}) we get eventually
$\mathcal G_{k, \ell} = -\alpha_{k-1} \rho_k^L \dots \rho_{\ell-1}^L\alpha_\ell^\dagger$. 
\QED
\subsection{Proof of Proposition \ref{new1}}
First, let us rephrase the computation of $\varphi_1^R$. We look for  $\xi = [\xi_1 \cdots \xi_p]\in \mathbb M_{N,p}$ ,  
{\it unitary} in the sense that $\ll \xi, \xi\gg= I_p$, 
"orthogonal" to $\varepsilon$ in the sense that  $\ll \xi, \varepsilon\gg= 0_p$ and 
such that the vectors of $\xi$ belong to $\spann\{\varepsilon, U\varepsilon\}$.
In a first step, let us see that the matrix $U\varepsilon - \varepsilon \gamma$ is orthogonal to $\varepsilon$ if and only if
\[\ll U\varepsilon, \varepsilon \gg\!\ =\!\  \ll \varepsilon\gamma, \varepsilon\gg\!\ =\!\ \gamma^\dagger\ll\varepsilon, \varepsilon\gg\]
hence $\gamma = \alpha^\dagger$. Let us now normalize this vector. The square of its "norm" is 
\[\ll U\varepsilon - \varepsilon  \alpha^\dagger, U\varepsilon - \varepsilon  \alpha^\dagger\gg\!\  = I_p - \alpha\alpha^\dagger\,,\]
so that, using the notation of Simon for the {\it defect} matrices
we claim that the matrix  
\ben
\label{defxi}\xi = \left(U\varepsilon - \varepsilon \alpha^\dagger\right)(\rho^R)^{-1}\een
satisfies all the requirements above. Of course, we demand that $\rho^R$ is invertible, but it is true in the generic case.

As in the scalar case, we define now an endomorphism $V$ unitary letting invariant the subspace orthogonal to 
$\spann\{\varepsilon, U\varepsilon\}$ and such that $V\varepsilon = U\varepsilon$. 
We know already from (\ref{defxi}) that $V\varepsilon = \xi \rho^R + \varepsilon \alpha^\dagger$. 
In the "basis" $(\varepsilon , \xi)$, we can say that the matrix of the restriction of $V$ will be $\Theta(\alpha)$, as defined in (\ref{deftheta0}),
 in the sense that if $w= \varepsilon a +  \xi b$ then $Vw =  \varepsilon a' +  \xi b'$ with
\[a'= \alpha^\dagger a + \rho^L b \ , \ b'= \rho^R a - \alpha b\,.\]

Now, the endomorphism 
\ben\label{defw}W = V^{-1}U\een is unitary and it fixes $\varepsilon$.  
In the basis obtained by orthonormalization of $\{\varepsilon, U \varepsilon, U^2\varepsilon, \dots\}$, the endomorphism $U$ has the block GGT matrix $\mathcal G^R(\alpha_0, \dots)$. In this basis $V(U, \varepsilon)$ has the matrix $\widetilde \Theta(\alpha_0)$. and by  (\ref{AGRGGT}) the restriction of $W(U, \varepsilon)$ to $\varepsilon^\perp$  has the matrix $\mathcal G^R(\alpha_1, \dots)$.
\
\QED

\subsection{Proof of Proposition \ref{propsursoi}}
\label{supropsursoi}
The first assertion  is a straightforward consequence of the invariance of the Haar measure. 

To prove the second assertion,  
we will follow some notation of Collins doctoral dissertation
 (\cite{CollinsHAL}) Section 4.2.
Let $\pi$ 
 be the canonical projection
: $\mathbb M_{N,N}  \rightarrow \mathbb M_{p, N}$. 
The set $\pi(\mathbb U(N))$ is a real sub-manifold of  $\mathbb M_{p,N}$ of dimension $p(2N-p)$, characterized by
\ben
\label{ch}
\pi(\mathbb U(N)) = \{V \in \mathbb M_{p, N} : V V^\dagger = I_p\}\,.\een 
The pushforward $\widehat\pi$ of the Haar measure on $\mathbb U(N)$ by $\pi$
is invariant under the natural action at left and right of  $\mathbb U(p)$ and $\mathbb U(N)$, respectively.  
Since the action of $\mathbb U(p) \times \mathbb U(N)$ on $\pi(\mathbb U(N))$ is transitive, this measure is the only normalised invariant one. 

For $M \in \mathbb M_{N,N}$ set 
\ben
\label{pipistar}
h(M) = \pi(M)\pi(M)^\dagger\een and   \ben
\label{vrac}v(M) = h(M)^{-1/2}\pi(M)\een (when $h(M) > 0$), 
where the square root is  taken in the functional calculus sense.

From (\ref{ch}), $v(M) \in \pi(\mathbb U(N))$. 
Let us now provide  $\mathbb M_{N,N}$ with the distribution $\tt Gin$$(N)$
 and let us denote by $\widetilde \pi$  the pushforward of this measure by  $v$. Let us prove that
\[\widetilde \pi = \widehat\pi\,.\]
Ii is enough to show that   $\widetilde \pi$ is left and right invariant.
Since $\pi(MU) = \pi(M)U$ for $U \in \mathbb U(N)$ and then $h(MU) = h(M)$,
and since the Gaussian distribution is invariant by $U$, the right invariance by $\mathbb U(N)$ is obvious.
Let us consider the left invariance.
If $U \in \mathbb U(p)$ and if $\widetilde U$ is defined by 
 \[ \widetilde U = \begin{pmatrix} U & 0_{p, N-p}\\
 0_{N-p, p} & I_{N-p}\end{pmatrix}\,,\]
 then $U\pi(M) = \pi(\widetilde U M)$ , $h(M) = U^\dagger h(\widetilde U M) U$ , $h(M)^{-1/2} = U^\dagger (h(\widetilde U M))^{-1/2} U$
and eventually
$U v(M) = v(\widetilde U M)$. The invariance of the  Gaussian distribution by $\widetilde U$ ends the job.
(Let us notice that  (\ref{pipistar}) is precisely the  relation (\ref{H=}). ) \QED
\subsection{Proof of Theorem \ref{alphajlaw}}
There are two approaches in the scalar case, that of  \cite{Killip1} and that of   \cite{simon2006cmf} section 11. We follow the method of proof of Theorem 11.1 in  \cite{simon2006cmf} 
 and extend it to the matricial case. The only difficulty comes from the noncommutativity.
%
%

Let $\varepsilon =[\varepsilon_1\cdots\varepsilon_p]$ an orthonormal system and (see  (\ref{relalpha0U})
\ben
\label{defalpha}
\alpha (U, \varepsilon) =\!\ \ll U\varepsilon, \varepsilon\gg\,.
\een
 If $U$ is Haar distributed, we have to find the distribution of $\alpha$ and to check that conditionally upon $\alpha$, the matrix $W$ is Haar distributed on $\mathbb U(N-p)$. Actually, $\alpha$ is nothing else than the upper left corner of size $p$ of $U$ and its distribution is known from Collins to be the $\tt Cor  (N,p)$ one. To prove the remaining part, let us see how the different quantities depend on $U$. To be clear, let us write  $\alpha(U , \varepsilon)$ for $\alpha$ and the same for $\xi,V$ and $W$ as defined in the proof of Proposition 2.8..
To characterize the Haar distribution, we use the criterion of invariance by left multiplication by an unitary matrix. We have to prove that if 
$\Gamma$ is a fixed unitary matrix letting $\varepsilon$ invariant and if $F$ is a Borel function 
\ben
\label{tobeproved}
\mathbb E F( W(U,\varepsilon)\Gamma, \alpha(U,\varepsilon)) = \mathbb E F (W(U,\varepsilon), \alpha(U,\varepsilon))\,.
\een
We have
\[ W(U,\varepsilon)\Gamma =  V(U,\varepsilon)^{-1}U\Gamma\]
which leads to consider the unitary endomorphism $U' = U\Gamma$  and to check successively
\[\alpha(U', \varepsilon) = \alpha(U, \varepsilon) \ , \  \xi(U', \varepsilon)= \xi(U, \varepsilon)\,,\]
(see (\ref{defalpha}) and (\ref{defxi})). 
It should be clear that $V(U', \varepsilon) = V(U, \varepsilon)$ since they coincide on  $\spann \{\varepsilon , \xi\}$ (see (\ref{deftheta0}) and leave invariant its orthogonal subspace. We have then
\[W(U, \varepsilon)\Gamma = V(U', \varepsilon)^{-1} U' = W(U', \varepsilon)\]
and then
\[\mathbb E F( W(U,\varepsilon)\Gamma, \alpha(U,\varepsilon)) =  \mathbb E F( W(U',\varepsilon), \alpha(U',\varepsilon))\,.\]Since $U'$ and $U$ have the same distribution, we have checked (\ref{tobeproved})  and the proof of the first step of the iteration is complete. 
Now we have to consider a matrix Haar distributed in $\mathbb U(N-p)$. It's the same reasoning. We stop the recursion at $j=J$.
\QED
\subsection{Proof of Corollary \ref{corharry}}
\label{proofharry}
With the notations of  Sec. \ref{suuu}, $\mu_{\w,p}$ is the spectral measure of the pair $(U, \mathcal H_p)$ so that, for $p$ and $n_0$ fixed, 
\begin{equation}
\label{canonical}
(U^n)_{i,j=1}^p = m_n (\mu_{\w,p})\ , \ \  1 \leq n \leq n_0\,.\end{equation}
%
%
Now, the $n^{\mbox{th}}$ matrix Verblunsky coefficient  $\alpha_n\sn$ of $\mu_{\w,p}$ is  also the $(n+1)$th matrix canonical moment 
 of $\mu_{\w,p}$ (as soon as  the dimension $N$ of the ambient space is greater than $ (n+2)p$) and   there is  a Taylor expansion of the moments 
$m_1, \dots, m_{n_0})$ 
 in function of the canonical moments $\alpha_0, \dots , \alpha_{n_0-1}$ (see 
 \cite{118218FGAR} Lemma 3.2). 
 Then it is enough to use the  so-called delta-method to deduce the statement of Corollary \ref{corharry} from the second statement of Proposition \ref{CLT} . Notice that for $n_0 =1$ the result goes back to \cite{PetzReffy}.
\QED
\subsection{Proof of Theorem \ref{rateptheo}}
\label{suproLDPVer}
There is actually two possible proofs. The first one (that is presented here), is short and use directly the Verblunsky coefficients. The second one is quite longer and does not use the Verblunsky coefficients but directly the representation on the eigenvalues of the matrix measure. This second proof is much more general as it may be applied to a general sequence of matrix-valued random measures and is useful to obtain general sum rules. This point of view will be developed in the forthcoming paper \cite{GaNaRo}.\\
\\
First of all, invoking Dawson-Gärtner's theorem on projective limits (\cite{demboZ98} Th.4.6.1),  we get the LDP for the random matrix measure at scale $N$ with good rate function
\[I(\mu) = -\sum_{j=0}^\infty \log \det (I-\alpha_j^\dagger\alpha_j)\,.\]
To conclude, we use the matricial Verblunsky form of the Szeg\H{o} theorem (see originally \cite{delsarte}, \cite{simon05} Theorem 2.13.5 and more recently \cite{derevy}):
\[I(\mu) = -\frac{1}{2\pi} \int_{-\pi}^\pi \log\det W(\theta) d\theta \] 
if $d\mu(\theta) = W(\theta)\frac{d\theta}{2\pi} +  d\mu_s(\theta)$ (Lebesgue decomposition). Notice that the last expression early appears in \cite{widom1975limit} in the asymptotic expansion of determinant of block Toeplitz matrices (see also \cite{dym1981abstract} for related results on more general block operators). 
\QED

\subsection{Proof of Corollary \ref{rem}}
Starting from
\[\nu = \nu_a  + \nu_s \ , \  \nu_a=\nu'_a dz\]
we see that $\nu_a \ll dz$  yields $\nu_a \lll I_p dz$. 
If $\nu'_a > 0$ a.e. then, by theorem 5.5 in \cite{robertson} we have $I_p dz \lll \nu_a$, hence $I_p dz \lll \nu$ and the Radon-Nikodym derivative of $I_p dz$ with respect  to $\nu$ is $(\nu'_a)^{-1}$ and (\ref{ratep0}) is valid.

Conversely, if $I_p dz \lll \nu$, then there exists a finite measure  $\gamma$ on $\mathbb T$ such that $I_p dz \ll \gamma$, $\nu \ll \gamma$ and range $(I_p dz/d\gamma) \subset$ range $(\nu'_\gamma)$. But  $I_p dz \ll \gamma$ implies $dz \ll \gamma$, so that $dz = g(z)d\gamma$. Since $I_p dz/ d\nu = g(z) (\nu'_\gamma(z))^\sharp$ the finiteness of the integral in (\ref{ratep0}) has two consequences:
\begin{itemize}
\item 
$g(z) >0$  for a.e. $z$ and then
\[d\gamma =  (g(z))^{-1} 1_{g(z) > 0} dz + d\gamma_s \ \ \ \gamma_s \perp dz\]
\item
$(\nu'_\gamma)^\sharp(z)\not=0$  for a.e. $z$. From the definition of the pseudo-inverse, this last requirement needs  $\nu'_\gamma\not=0$  for a.e. $z$.
\end{itemize}
 This yields
\[d\nu = \nu'_\gamma d\gamma = \nu'_\gamma(z) (g(z))^{-1} 1_{g(z) > 0} dz + \nu'_\gamma(z) d\gamma_s\]
and then $\nu'_a (z) = \nu'_\gamma(z) (g(z))^{-1}$ for a.e. $z$ and (\ref{ratep}) is valid. \QED
\subsection{Proof of Theorem \ref{theoarem}}
\label{suprothediminf}
As we noticed in (\ref{proj})  the structure of spectral measures is projective. We may apply the Dawson-Gärtner theorem (\cite{demboZ98} Th. 4.6.1) and we get the rate function 
\begin{equation}
\label{fdet}\sup_k \int_\mathbb T - \log \det {\mu'_a}^k(z) dz\,.\end{equation} 
Then we use the following lemma whose proof is slightly postponed.\QED

\begin{lem}
\label{lem7}
$\;$
In the settings of the theorem the following statements hold true:
\begin{enumerate}
\item 
The sequence $\int_\mathbb T - \log \det {\mu'_a}^k(z) dz $ is increasing in $k$.
\item We have
\begin{equation}\label{Hajek} \lim_k  \det {\mu'_a}^k(z) = \det \mu'_a(z)\,.\end{equation}
\end{enumerate} 
\end{lem}
The first statement entails that  the supremum in (\ref{fdet}) is actually an increasing limit.
The second statement gives a limit for the integrand. But, in
 general, it is not possible to commute limit and integral in (\ref{iinfini}).  Assumption (\ref{C}) ensures a dominated convergence.
The equality $\tr \log = \log \det$ is classical (see \cite{derevy}).\\
\\
{\bf Proof of Lemma \ref{lem7} }
\begin{enumerate}
\item For fixed $k$, the   Hermitian non-negative matrix ${\mu'_a}^k(z)$ admits a Cholesky decomposition
\begin{equation}
\label{chol}{\mu'_a}^k(z) = L_k(z)L_k(z)^\dagger\,,\end{equation}
and it is straightforward to see that the $(k-1)$-section  of $L_k(z)$ is $L_{k-1}(z)$, so that we have an infinite Cholesky matrix $L(z)$ whose  generic  entry
 will be denoted by $\ell_{i,j}(z)$. From (\ref{chol}), we have the relation
\begin{equation}
\label{bartlett}
\det {\mu'_a}^k(z) = \det{\mu'_a}^{k-1}(z) |\ell_{kk}(z)|^2\,.
\end{equation}
Taking logarithm and integrating the last relation in $z$ we get
\begin{equation}
\label{justelog}\int_\mathbb T -\log \det  {\mu'_a}^k(z) dz + \int_\mathbb T \log \det  {\mu'_a}^{k-1}(z)  = \int_\mathbb T -\log |\ell_{kk}(z)|^2 dz\end{equation}
which, by Jensen's inequality gives 
\begin{equation}
\label{jensen} 
\int -\log \det  {\mu'_a}^k(z) dz + \int \log \det  {\mu'_a}^{k-1}(z) \geq -\log \int |\ell_{kk}(z)|^2 dz\end{equation}
From (\ref{chol}) we have  also
\begin{equation}
\label{chol2}
|\ell_{kk}(z)|^2 = 
\left(\mu'_a(z)\right)_{kk} -\sum_{j=0}^{k-1}|\ell_{kj}(z)|^2 \leq \left(\mu'_a(z)\right)_{kk}\,.
\end{equation}
Now, $\mu$ is a matrix probability measure and then $\int_{\mathbb T} \mu'_a(z) dz \leq I$, which implies 
$\int \left(\mu'_a(z)\right)_{kk} dz \leq 1$, and by integration in (\ref{chol2}) 
\begin{equation}
 \int_\mathbb T |\ell_{kk}(z)|^2 dz \leq 
 1\,.
\end{equation}
Plugging in (\ref{jensen}) ends the proof of the first part of  Lemma \ref{lem7}.

\item 
This follows directly from the Cholesky decomposition in (\ref{chol}).
\end{enumerate}
\QED
{\bf Acknowledgment}\\
The authors would like to warmly thank the anonymous reviewers
for their valuable comments and suggestions to improve the
quality of the paper.

\bibliographystyle{plain}
\bibliography{blockbib}

\begin{thebibliography}{10}

\bibitem{Arlinski}
Y.~Arlinskii.
\newblock Conservative discrete time-invariant systems and block operator {CMV}
  matrices.
\newblock {\em Methods Funct. Anal. Topology}, 15(3):201--236, 2009.

\bibitem{bai1}
Z.~D. Bai, B.~Q. Miao, and G.~M. Pan.
\newblock On asymptotics of eigenvectors of large sample covariance matrix.
\newblock {\em Ann. Probab.}, 35(4):1532--1572, 2007.

\bibitem{Changetal}
F.C. Chang, J.H.B. Kempermann, and W.J. Studden.
\newblock A normal limit theorem for moment sequences.
\newblock {\em Ann. Probab.}, 21(3):1295--1309, 1993.

\bibitem{CollinsHAL}
B.~Collins.
\newblock {\em Intégrales matricielles et Probabilités non-commutatives}.
\newblock PhD thesis, Paris 6, 2004.
\newblock {\it Available at}
  \href{http://tel.archives-ouvertes.fr/tel-00004306/fr/}{tel-00004306/fr}{(in
  english)}.

\bibitem{collins2005product}
B.~Collins.
\newblock {Product of random projections, {J}acobi ensembles and universality
  problems arising from free probability}.
\newblock {\em Probab. Theory Related Fields}, 133(3):315--344, 2005.

\bibitem{damanik2008analytic}
D.~Damanik, A.~Pushnitski, and B.~Simon.
\newblock {The analytic theory of matrix orthogonal polynomials}.
\newblock {\em Surv. Approx. Theory}, 4:1--85, 2008.

\bibitem{delsarte}
P.~Delsarte, Y.V. Genin, and Y.G. Kamp.
\newblock Orthogonal polynomial matrices on the unit circle.
\newblock {\em IEEE Trans. Circuits and Systems}, pages 149--160, 1978.

\bibitem{demboZ98}
A.~Dembo and O.~Zeitouni.
\newblock {\em Large Deviations Techniques and Applications}.
\newblock Springer, 1998.

\bibitem{derevy}
M.~Derevyagin, O.~Holtz, S.~Khrushchev, and M.~Tyaglov.
\newblock Szeg{\H o}'s theorem for matrix orthogonal polynomials.
\newblock {\em J. Approx. Theory}, 164(9):1238--1261, 2012.

\bibitem{DeGa}
H.~Dette and F.~Gamboa.
\newblock Asymptotic properties of the algebraic moment range process.
\newblock {\em Acta Math. Hungar.}, 116(3):247--264, 2007.

\bibitem{denag09}
H~Dette and J.~Nagel.
\newblock Matrix measures, random moments, and {G}aussian ensembles.
\newblock {\em J. Theoret. Probab.}, 25(1):25--49, 2012.

\bibitem{DeSt97}
H.~Dette and W.J. Studden.
\newblock {\em The theory of canonical moments with applications in statistics,
  probability, and analysis}.
\newblock Wiley Series in Probability and Statistics,, 1997.

\bibitem{destu02}
H.~Dette and W.J. Studden.
\newblock Matrix measures, moment spaces and {F}avard's theorem for the
  interval [0,1] and [0,$\infty$).
\newblock {\em Linear Algebra Appl.}, 345:169--193, 2002.

\bibitem{destu03}
H.~Dette and W.J. Studden.
\newblock Quadrature formulas for matrix measures---a geometric approach.
\newblock {\em Linear Algebra Appl.}, 364:33--64, 2003.

\bibitem{dewag10}
H.~Dette and J.~Wagener.
\newblock Matrix measures on the unit circle, moment spaces, orthogonal
  polynomials and the {G}eronimus relations.
\newblock {\em Linear Algebra Appl.}, 432:1609--1626, 2010.

\bibitem{DuSch}
N.~Dunford and J.T. Schwartz.
\newblock {\em Linear operators. {P}art {II}}.
\newblock Wiley Classics Library. John Wiley \& Sons Inc., New York, 1988.
\newblock Spectral theory. Selfadjoint operators in Hilbert space, With the
  assistance of William G. Bade and Robert G. Bartle, Reprint of the 1963
  original, A Wiley-Interscience Publication.

\bibitem{dym1981abstract}
H.~Dym and S.~Ta'assan.
\newblock An abstract version of a limit theorem of \mbox{Szeg{\"o}}.
\newblock {\em Journal of Functional Analysis}, 43(3):294--312, 1981.

\bibitem{Gamboa}
F.~Gamboa and L.-V. Lozada-Chang.
\newblock Large deviations for random power moment problem.
\newblock {\em Ann. Probab.}, 32(3B):2819--2837, 2004.

\bibitem{GaNaRo}
F.~Gamboa, J.~Nagel, and A.~Rouault.
\newblock From large deviations to sum rules, 2012.
\newblock {\it Work in progress}.

\bibitem{118218FGAR}
F.~Gamboa, J.~Nagel, A.~Rouault, and J.~Wagener.
\newblock Large deviations for random matricial moment problems.
\newblock {\em J. Multivariate Anal.}, 106:17--35, 2012.

\bibitem{gamboacanonical}
F.~Gamboa and A.~Rouault.
\newblock {Canonical moments and random spectral measures}.
\newblock {\em J. Theoret. Probab.}, 23:1015--1038, 2010.

\bibitem{FGAR}
F.~Gamboa and A.~Rouault.
\newblock Large deviations for random spectral measures and sum rules.
\newblock {\em Appl. Math. Res. Express}, pages 281--307, 2011.

\bibitem{Gup}
J.C. Gupta.
\newblock Completely monotone multisequences, symmetric probabilities and a
  normal limit theorem.
\newblock {\em Proc. Indian Acad. Sci. Math. Sci.}, 110:415--430, 2000.

\bibitem{KarStu}
S.~Karlin and W.~J. Studden.
\newblock {\em Tchebycheff systems: {W}ith applications in analysis and
  statistics}.
\newblock Pure and Applied Mathematics, Vol. XV. Interscience Publishers John
  Wiley \& Sons, New York-London-Sydney, 1966.

\bibitem{Killip1}
R.~Killip and I.~Nenciu.
\newblock Matrix models for circular ensembles.
\newblock {\em Int. Math. Res. Not.}, 50:2665--2701, 2004.

\bibitem{KreNud}
M.~G. Krein and A.~A. Nudel{$'$}man.
\newblock {\em The {M}arkov moment problem and extremal problems}.
\newblock American Mathematical Society, Providence, R.I., 1977.
\newblock Ideas and problems of P. L. Chebyshev and A. A. Markov and their
  further development, Translated from the Russian by D. Louvish, Translations
  of Mathematical Monographs, Vol. 50.

\bibitem{krishna09}
M.~Krishnapur.
\newblock From random matrices to random analytic functions.
\newblock {\em Ann. Probab.}, 37:314--346, 2009.

\bibitem{LozEJP}
L.V. Lozada-Chang.
\newblock Large deviations on moment spaces.
\newblock {\em Electronic J. Probab.}, 10:662--690, 2005.

\bibitem{mandrekar1}
V.~Mandrekar and H.~Salehi.
\newblock On singularity and {L}ebesgue type decomposition for operator-valued
  measures.
\newblock {\em J. Multivariate Anal.}, 1(2):167--185, 1971.

\bibitem{Ner1}
Y.A. Neretin.
\newblock Hua-type integrals over unitary groups and over projective limits of
  unitary groups.
\newblock {\em Duke Math. J.}, 114(2), 2002.

\bibitem{PetzReffy}
D.~Petz and J.~R{\'e}ffy.
\newblock On asymptotics of large {H}aar distributed unitary matrices.
\newblock {\em Period. Math. Hungar.}, 49(1):103--117, 2004.

\bibitem{robertson}
J.B. Robertson and M.~Rosenberg.
\newblock {The decomposition of matrix-valued measures}.
\newblock {\em Michigan Math. J}, 15:353--368, 1968.

\bibitem{Al1}
A.~Rouault.
\newblock Asymptotic behavior of random determinants in the {L}aguerre, {G}ram
  and {J}acobi ensembles.
\newblock {\em ALEA Lat. Am. J. Proba Math. Stat.}, 3:181--230, 2007.

\bibitem{simon05}
B.~Simon.
\newblock {\em Orthogonal polynomials on the unit circle. Part 1: Classical
  theory and Part 2: Spectral Theory.}
\newblock {Colloquium Publications. American Mathematical Society, Providence,
  RI: American Mathematical Society (AMS)}, 2005.

\bibitem{simon2006cmf}
B.~Simon.
\newblock C{MV} matrices: {F}ive years after.
\newblock {\em J. Comput. Appl. Math.}, 208:120--154, 2007.

\bibitem{widom1975limit}
H.~Widom.
\newblock On the limit of block toeplitz determinants.
\newblock {\em Proceedings of the American Mathematical Society}, pages
  167--173, 1975.

\end{thebibliography}
\end{document}